\documentclass[a4paper,reqno]{amsart}
\pdfoutput=1

\usepackage{amsmath,amsfonts,amssymb,amsthm,stmaryrd,bbm,graphicx,mathtools,enumerate}
\usepackage{mathrsfs}
\usepackage[T1]{fontenc} 
\usepackage[latin1]{inputenc}
\usepackage[english]{babel}
\usepackage[tracking,spacing,kerning,babel]{microtype}
\usepackage{wasysym} 
\usepackage{color}
\usepackage{xcolor}
    	\usepackage{framed} 
\usepackage{wrapfig} 

\usepackage[final]{hyperref}   
\hypersetup{
    linktoc=page,
    linkcolor=red,          
    citecolor=blue,        
    filecolor=blue,      
    urlcolor=cyan,
   colorlinks=true           
}

\usepackage{lipsum} 

\usepackage{minitoc}
\usepackage{multicol}

\newtheoremstyle{mine}
{\baselineskip}
{\baselineskip}
{\itshape}
{
}
{\bfseries}
{.}
{.5em}
{#1 #2\ifx#3\relax\else~(#3)\fi}

\theoremstyle{mine}

\newtheorem{theorem}{Theorem}
\numberwithin{theorem}{section}
\newtheorem{corollary}[theorem]{Corollary}
\newtheorem{proposition}[theorem]{Proposition}
\newtheorem{lemma}[theorem]{Lemma}

\newtheorem{definition}[theorem]{Definition}
 
\numberwithin{equation}{section}
\newtheorem{remark}[theorem]{Remark}

\colorlet{shadecolor}{blue!10}

\renewcommand{\d}[1]{\ensuremath{\operatorname{d}\!{#1}}}

\renewcommand{\epsilon}{\varepsilon}

\newcommand{\R}{\mathbb{R}}

\newcommand{\N}{\mathbb{N}}

\renewcommand{\d}{\mathbf d}

\def\calL{\mathcal{L}}

\def\calN{\mathcal{N}}

\def\P{\mathbb{P}} 
\def\E{\mathbb{E}} 
\def\md{\mid}

\def\Bb#1#2{{\def\md{\bigm| }#1\bigl[#2\bigr]}}

\def\Eb{\Bb\E}

\def\<#1{\langle #1\rangle}

\def\bi{\begin{itemize}}  
\def\ei{\end{itemize}}
\def\bnum{\begin{enumerate}} 
\def\enum{\end{enumerate}}
\def\ni{\noindent}

\def\PPP{\mathrm{PPP}}

\newcommand{\ind}[1]{\mathbf{1}_{\left\{#1\right\}}}

\title
[A new proof of Liggett's theorem]
{
A new proof of Liggett's theorem for non-interacting Brownian motions
}

\author{Xinxin Chen, Christophe Garban,  Atul Shekhar}

\address
{Universit\'e Claude Bernard Lyon 1, CNRS UMR 5208, Institut Camille Jordan, 69622 Villeurbanne, France}
\email{chen@math.univ-lyon1.fr; shekhar@math.univ-lyon1.fr}

\address
{Universit\'e Claude Bernard Lyon 1, CNRS UMR 5208, Institut Camille Jordan, 69622 Villeurbanne, France \, and Institut Universitaire de France (IUF)}
\email{garban@math.univ-lyon1.fr}

\begin{document}

\maketitle

\begin{abstract}
In this note, we give a new proof of Liggett's theorem on the invariant measures of independent particle systems from \cite{liggett78} in the particular case of independent drifted Brownian motions. This particular case has received a lot of  attention recently due to its applications for the analysis of the local extrema of discrete Gaussian free field. 
The novelty of our proof is that it identifies directly the expected Poisson Point Process with exponential intensity without  relying on the Choquet-Deny convolution equation $\mu * P=\mu$ (\cite{ChoquetDeny,Deny}).
\end{abstract}

\section{Introduction}

\subsection{Context and main result.}
In his seminal paper \cite{liggett78},  Liggett gave a full characterisation of point processes which are left invariant under  non-interacting particle systems. His setting is very general (we refer to \cite{liggett78} as well as to \cite[Chapter 9]{biskupBook} for a very clear exposition of Liggett's theory). Let us shortly describe what his general result is. Let $\mathscr{H}$ be a locally compact, separable Hausdorff space and fix $P$ to be a Markov kernel on $\mathscr{H}$. Liggett considers the following discrete-time dynamics  on (locally finite) point processes $\theta$ on $\mathscr{H}$. Starting from an initial point process $\theta_0=\theta$, particles evolve independently of each other as  discrete-time Markov chains on $\mathscr{H}$ with kernel $P$. Under very mild conditions, Liggett proves that the invariant measures of this non-interacting particle system $n\mapsto \theta_n$ are classified by the Radon measures $\mu$ on $\mathscr{H}$ which satisfy the celebrated convolution equation of Choquet-Deny 
(\cite{ChoquetDeny,Deny})
\[
\mu*P=\mu\,.
\]
More precisely he shows that all the invariant measures are obtained as superpositions of Poisson Point Processes whose intensity measures are given by the solutions to Choquet-Deny convolution equation. 
Liggett's theorem in the particular case where $\mathscr{H}=\R$ and where particles evolve as independent drifted Brownian motions\footnote{It is more natural in this particular case to work in continuous time} has known an important revival over the last years especially since the work by Biskup and Louidor \cite{BL16}. In this work, the authors relied on Liggett's theorem to show, via a  beautiful ``Dysonization procedure'',  that local extrema of a Discrete Gaussian Free field are asymptotically distributed as a shifted Poisson Point Process with intensity $e^{-\lambda x} dx$.  See also \cite[Chapter 9]{biskupBook} for a very nice account on the characterization by Dysonization as well as \cite{zeitouni2017} where such a Dysonization procedure is also used. 
Another important use of Liggett's theorem, also in the case $\mathscr{H}=\R$ case, can be found in the works by Ruzmaikina-Aizenman and Arguin-Aizenman \cite{AR05,AA09} where links with spin glasses are highlighted.  

In this paper, we focus on this case of drifted Brownian motions on $\R$ which is thus relevant to \cite{AR05,AA09,BL16,zeitouni2017}. Liggett's non-interacting particle system corresponds in this case to the following process. 
Let $\theta$ be a locally finite point process on $\R$. (See the space $\calN$ and its topology below in Section \ref{s.prel}). Let us write 
\[
\theta = \sum_{i\in I}\delta_{x_i},
\]
with some finite or countable index set $I$. To each atom $x_i$ of $\theta$, we attach an independent Brownian motion with drift $-\lambda$ which is denoted as $(B_t^i-\lambda t)_{t\ge0}$. (N.B. the minus sign here is just a convention, the drift $\lambda$ may take any value in $\R$). At time $t>0$, we get the following point process $\theta_t$
\[
\theta_t=\sum_{i\in I}\delta_{x_i+B_t^{i}-\lambda t}.
\]
As pointed out already in \cite{liggett78} the fact that $\theta_0$ is locally finite does not imply a priori that $\theta_t$ still is. (It is easy to build examples for which a {\em coming down from $\infty$} happens at time $0+$). 
Because of these possible explosions, our notion of invariant measure is not in terms of Feller processes say, but rather following \cite{liggett78} in the following sense: a point process $\theta\sim \pi$ (where $\pi$ is a probability measure on the space $(\calN,d)$, see Section \ref{s.prel}) will be said to be an invariant measure (or a fixed point) for Brownian motion with drift $-\lambda$ if for all $t>0$, 
\[
\theta_t \overset{d}{=}\theta\,.
\]
See also our companion paper \cite[Section 2.2]{CGS-BBM} for a detailed discussion on this notion of invariant measure. 
We may now state Liggett's theorem in this particular setting. 
\begin{theorem}[Theorem $1.2$ of \cite{liggett78}]\label{main-thm}
For any $\lambda\in \R$,  a point process $\theta$ is invariant under Brownian Motions with drift $-\lambda$ if and only if $\theta$ is distributed as a Poisson point process with a random intensity measure\footnote{These processes are often called Cox processes} of the form 
\begin{equation}\label{int-form}
(Z_\infty e^{-2\lambda x}+ Y_\infty)dx,
\end{equation}
where $Z_\infty$ and $Y_\infty$ are some non-negative (possibly correlated) random variables. 
\end{theorem}

Liggett's original proof of the above Theorem is based on Choquet-Deny Theorem \cite{Deny}, \cite{ChoquetDeny} on convolution equations. See  Section \ref{recap} below for a brief summary of Liggett's proof. Our new proof avoids the use of Choquet-Deny Theorem and in some sense computes the solutions of Deny-Choquet equation along the way. Our method would apply for more general Markov Kernels on $\R^d$ than the one corresponding to drifted Brownian motions on $\R$, but for simplicity we shall stick to this case which seems to be the most relevant one in the literature. 

More than the result itself (which is a particular case of \cite{liggett78}), we believe that the main interest of this note lies in the strategy of identification of the invariant measures. For example in  \cite{CGS-BBM}, we were unable to adapt Liggett's proof scheme in order to characterize the invariant measures of branching Brownian motion with critical drift. See the discussion in Section 1.5 in \cite{CGS-BBM}. We therefore had to look for a different approach which is presented here in the simpler setting of drifted Brownian motions without branching. We believe this new way of identifying fixed points may be applied to other natural settings.


\section{Short preliminaries}\label{s.prel}

\ni
\textbf{State space.}
Let $\calN$ be the space of integer valued measures on $\R$ which are locally finite. 
This space is naturally equipped with the vague topology, see \cite{kallenberg2006}.
\begin{remark}\label{}
Note that the weak topology is not appropriate for the type of processes considered in this paper. This is due to the following reason: recall $\theta_n \overset{w}\longrightarrow \theta$ if and only if for any continuous bounded $f\in C_b(\R)$,
$\theta_n(f) \to \theta(f)$. But the processes we consider will in general have a diverging mass near $-\infty$, as such they will not integrate, say the continuous function $f\equiv 1$. The vague topology is more indulgent and corresponds instead to $\theta_n\overset{v}\longrightarrow \theta$ if and only if for any $f\in C_c(\R)$, $\theta_n(f) \to \theta(f)$.  
\end{remark}

The vague topology on $\calN$ is metrizable and one can define a metric $d=d_{\calN}$ on $\calN$ such that the space $(\calN,d)$ is Polish (see Theorem A2.3 in \cite{kallenberg2006}). As such we may now consider probability measures on $\calN$ in the usual way. 

%
%
\medskip

\ni
\textbf{Cox processes.}
A useful class of random variables in $(\calN,d)$ are the so-called {\em Cox processes.} They correspond  to the random point processes $\theta\in \calN$ which are defined via the following two steps procedure: 
\bnum
\item Start by sampling a random $\sigma$-finite positive measure $\sigma$ on $(\R, \mathcal{B}(\R))$.
\item Then, given $\sigma$, sample $\theta$ to be a Poisson Point Process of intensity $\sigma$, i.e. 
$\theta\sim \mathrm{PPP}(\sigma)$.
\enum 
We recall that the Laplace transform of a Cox process is given by the following L\'evy-Khintchine formula. 
\begin{proposition}[Proposition $2.12$ of \cite{bovier-book} or p28 in \cite{kingman}]
Let $\theta$ be a Cox process with random intensity $\sigma$.  Then, for any non-negative measurable function $f$, 
\begin{equation}\label{LK}
\E[e^{-\<{f, \theta}}] = \E\biggl[\exp\biggl\{-\int_{-\infty}^{\infty}(1-e^{-f(x)})\sigma(\d x)\biggr\}\biggr],
\end{equation}
where
\[\langle f,\theta \rangle := \int_{-\infty}^\infty f(x)\theta(\d x).\]
\end{proposition}

\section{Liggett's proof and idea of the new proof}\label{recap}


\subsection{Summary of Liggett's proof.}
(See \cite{biskupBook} or our summary of Liggett's proof in \cite{CGS-BBM} which holds for the more general setting). 
To characterize the point processes $\theta$ invariant under Brownian motions with drift $-\lambda$, it suffices to check (see for example \cite{kallenberg-book}) that for all $f\in C_c^+(\R) $, 
\begin{equation}\label{eqn-for-eta}
\E[e^{-\langle f,\theta_t \rangle}] = \E[e^{-\langle f,\theta \rangle}].   
\end{equation}

Using some basic computations, one obtains
\begin{equation}\label{exact-comp}
\E[e^{-\langle f,\theta_t \rangle}] = \E\left[\exp\left(\int_{\R}\log \E[e^{-f(x+B_t-\lambda t)}]\theta(\d  x)\right)\right],
\end{equation}
where $\E[e^{-f(x+B_t-\lambda t)}]=1-o_t(1)$ uniformly for $x\in\R$.  So, 
\begin{equation}\label{replace-log}
   -\log \E[e^{-f(x+B_t-\lambda t)}] = (1+ o_t(1))\E[1-e^{-f(x+B_t-\lambda t)}] .
\end{equation}
Let $P^t(x,\d y)$ denote the transition probabilities of $B_t - \lambda t $. Then, by use of Fubini theorem, one gets essentially that 
\begin{align}\label{factorisation}
    \E[e^{-\langle f,\theta_t \rangle}] = &\E\left[\exp\left(-(1+o_t(1))\int_\R\int_\R (1-e^{-f(y)})P^t(x,\d y)\theta(\d x)\right)\right]\nonumber\\
    =&\mathbb{E}\biggl[ \exp\biggl\{ -(1+ o_t(1))\int_\R(1-e^{-f(y)})M_t(\d y) \biggr\}\biggr],
\end{align}
where $M_t(\d y) = \langle P^t(\cdot,\d y), \theta(\cdot) \rangle$. Letting $t\to\infty$ leads to
\[
\E[e^{-\langle f,\theta \rangle}]=\lim_{t\to\infty}\mathbb{E}\biggl[ \exp\biggl\{ -\int_\R(1-e^{-f(y)})M_t(\d y) \biggr\}\biggr].
\]
Let $h(y):=1- e^{-f(y)}$ with $f\in C_c^+(\R)$. Then it follows that for any $h\in C_c^+(\R)$ such that $0\leq h(y)<1$,  $\lim_{t\to \infty} \E[e^{-\<{h,M_t}}]$ exists. 
As this class of functions $h$ is sufficiently large, we deduce that  $M_t(\d y)$ converges in law to some locally finite random measure $M_\infty(\d y)$. Furthermore, 
\begin{equation}
    \E[e^{-\langle f,\theta \rangle}] = \mathbb{E}\biggl[ \exp\biggl\{ -\int(1-e^{-f(y)})M_\infty(\d y) \biggr\}\biggr]\,.
\end{equation}
This implies that $\theta$ must be a Cox process with random intensity measure given by $M_\infty(\d y)$. Using the constraint \eqref{eqn-for-eta} once more, it follows that $M_\infty P^t\overset{d}{=}M_\infty$. Furthermore, since the law of $B_t$ is supported on the whole $\R$, Corollary $3.8$ of Liggett \cite{liggett78} implies that $M_\infty P^t = M_\infty$ a.s. This is the well-known convolution equation of Choquet-Deny and the Choquet-Deny Theorem (see \cite{ChoquetDeny,Deny}) implies that $M_\infty(dy)$ must necessarily be of the form \eqref{int-form}. 

\begin{figure}[!htp]
\begin{center}
\includegraphics[height=9.5cm]{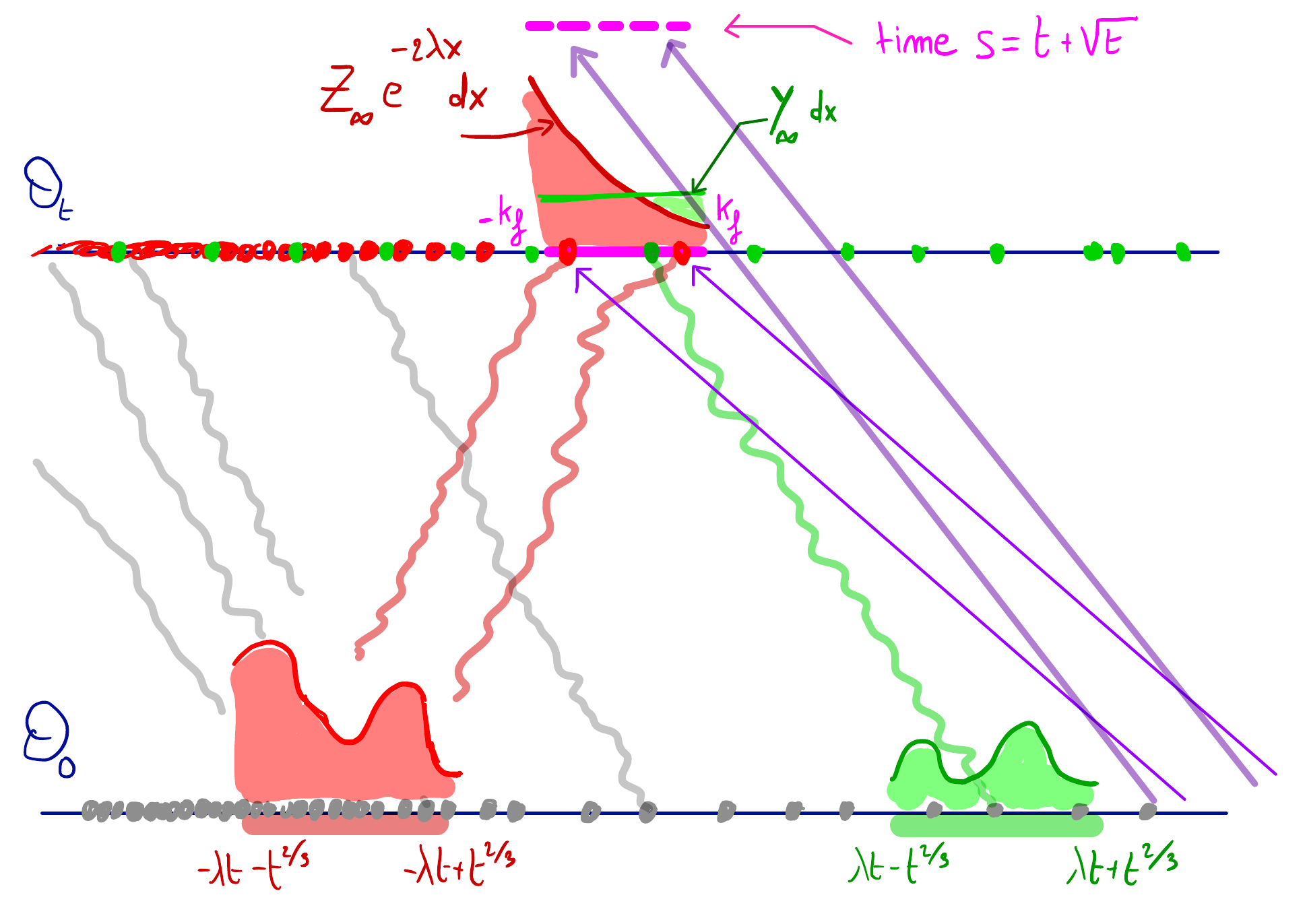}
\end{center}
\caption{An artistic view of the proof.}\label{f.sketch}
\end{figure}
\subsection{Idea of the new proof.}
To characterize the invariant point processes $\theta$, we still wish to recognize the structure of a Cox process via the Laplace transforms $\Eb{e^{-\<{f,\theta}}}$ for any fixed $f\in C_c^+(\R)$. Starting as in the previous analysis, for any $f\in C_c^+(\R)$ we have  as $t\to \infty$
\begin{align*}
    \E[e^{-\langle f,\theta_t \rangle}] &= \E\left[\exp\left(-(1+o_t(1))\int_\R\int_\R (1-e^{-f(y)})P^t(x,\d y)\theta(\d x)\right)\right]\nonumber \\
     &= \E\left[\exp\left(-(1+o_t(1))\int_\R \theta(\d x) \int_{[-K_f,K_f]} (1-e^{-f(y)})P^t(x,\d y) \right)\right]\,,\nonumber
\end{align*}
where $K_f$ is such that $\mathrm{supp}(f)\subset [-K_f,K_f]$. We are left with understanding the following random variable (possibly at large times $t$)  
\begin{align}\label{e.integral}
&\int_\R \theta(\d x) \int_{[-K_f,K_f]} (1-e^{-f(y)})P^t(x,\d y)\,.
\end{align}
Now the main observation in the proof is as follows: if all fixed points $\theta$ were indeed given (as we expect) by Cox processes with random intensity measure of the form~\eqref{int-form}, we would notice that in the above integral, as $t\to\infty$, with high probability, only points coming from very specific zones would contribute to~\eqref{e.integral}. Namely, for $t$ large 
\bi
\item Particles seen at time $t$ in the window $[-K_f,K_f]$ coming from the {\em flat} part $\PPP(Y_\infty dx)$ will arise with high probability from initial particles (i.e. at time $t=0$) in the far ahead region $[\lambda t-t^{2/3}, \lambda t + t^{2/3}]$. {\em N.B we may have chosen instead the window $[\lambda t -t^\alpha, \lambda t + t^\alpha]$ for any exponent $\alpha\in(\tfrac 12 , 1)$.}

\item While Particles seen at time $t$ in the window $[-K_f,K_f]$ coming from the {\em exponential} part $\PPP(Z_\infty e^{-2\lambda x} dx)$ will arise with high probability from initial particles ($t=0$) in the far backward region $[-\lambda t-t^{2/3}, -\lambda t + t^{2/3}]$. 
\ei 
This observation leads us to decompose the (random) integral~\eqref{e.integral} into the following three parts.
\begin{align*}
\mathbf{Z}_t(f):=&\int_{[-\lambda t-t^{2/3},-\lambda t+t^{2/3}]} \theta(\d x) \int_{[-K_f,K_f]} (1-e^{-f(y)})P^t(x,\d y)\\
\mathbf{Y}_t(f):=&\int_{[\lambda t-t^{2/3},\lambda t+t^{2/3}]} \theta(\d x) \int_{[-K_f,K_f]} (1-e^{-f(y)})P^t(x,\d y) \\
E_t(f):=&\int_{\R\setminus([-\lambda t-t^{2/3},-\lambda t+t^{2/3}]\cup [\lambda t-t^{2/3},\lambda t+t^{2/3}])}
\theta(\d x) \int_{[-K_f,K_f]} (1-e^{-f(y)})P^t(x,\d y)\,.
\end{align*}
The proof then proceeds by analyzing each of these three terms on which we now say a few words. See also Figure \ref{f.sketch} which illustrates what is the strategy of the proof. 
\bnum
\item The key point in the analysis of the part $\mathbf{Z}_t(f)$ will be to obtain estimates uniform in the starting point $x\in[-\lambda t-t^{2/3}, -\lambda t +t^{2/3}]$. Indeed, without having too much control on the way initial particles are spread over $[-\lambda t-t^{2/3}, -\lambda t +t^{2/3}]$, this will allow us to show that after a long time $t$ the few lucky particles among these which will reach a position in $[-K_f,K_f]$ will have ``equilibrated'' under an exponential distribution $1_{[-K_f,K_f]} e^{-2\lambda x}dx$. The only control we need on the initial spreading of particles in   $[-\lambda t-t^{2/3}, -\lambda t +t^{2/3}]$ is the fact that there are not too many particles there (with high probability at least). We use for this the fact that $\theta_t\overset{d}=\theta$. If there were too many particles at times $0$ in  $[-\lambda t-t^{2/3}, -\lambda t +t^{2/3}]$, we show that this would create a suspicious rise of the intensity at time $t$. 
\item The analysis of the part  $\mathbf{Y}_t(f)$ will follow the same idea  except that now particles come from the up-front $[\lambda t - t^{2/3}, \lambda t +t^{2/3}]$  (when $\lambda>0$, otherwise left and right need to be permuted). The only but important difference is that now the few particles which start at time $0$ in  $[\lambda t - t^{2/3}, \lambda t +t^{2/3}]$ and will reach $[-K_f,K_f]$ at time $t$ will now ``equilibrate'' according to the {\em flat} Lebesgue measure and are as such (asymptotically as $t\to \infty$)  responsible for the part $\PPP(Y_\infty dx)$.

\item To conclude the proof, it remains to show that it is very unlikely that particles initially in the sets 
\[
(-\infty, -\lambda t -t^{2/3}] \cup [-\lambda t+t^{2/3}, \lambda t -t^{2/3}] \cup [\lambda t + t^{2/3},\infty)  
\]
will reach at time $t$ the window $[-K_f,K_f]$. This is in fact the main step in the proof whose underlying idea is as follows: imagine that for arbitrary large times $t$, we find with probability bounded away from 0 particles issued say from the left interval $(-\infty, -\lambda t -t^{2/3}]$ . Then we claim that this would imply that at a well-chosen later time $s=t+\sqrt{t}$, there would be way too many particles in $[-K_f,K_f]$ which would contradict $\theta_s\overset{d}=\theta$. 
\enum

\section{A New proof of Theorem \ref{main-thm}.}

\subsection{Existence.}
We first check that Poisson point processes with random intensity measure $(Z_\infty e^{-2\lambda x}+ Y_\infty)dx$ are indeed invariant. This is a standard computation, and we provide it here for the sake of completeness. Let $\theta \overset{d}{=} \PPP((Z_\infty e^{-2\lambda x}+ Y_\infty)dx)$. Using the characterisation property of Laplace transforms (see for example Proposition $2.10$ of \cite{bovier-book}), it suffices to check that for all $f \in C_c^+(\R)$, 
\begin{equation}\label{lap-eq}
\E[e^{-\<{f, \theta_t}}] = \E[e^{-\<{f, \theta}}].
\end{equation} 
By \eqref{LK}, 
\[\E[e^{-\<{f, \theta}}] = \E\biggl[\exp\biggl\{-\int_{-\infty}^{\infty}(1-e^{-f(x)})(Z_\infty e^{-2\lambda x}+ Y_\infty)\d x\biggr\}\biggr].\]
On the other hand,  
\begin{align}\label{comp-1}
\nonumber\E[e^{-\<{f,\theta_t}}]=&\E\left[\exp\left\{\int_{-\infty}^\infty  \log(\E[e^{-f(x+B_t-\lambda t)}])    \theta(\d x)\right\}\right]\\
=&\E[\exp\{-\int_{-\infty}^\infty \E[1-e^{-f(x+B_t-\lambda t)}](Z_\infty e^{-2\lambda x}+ Y_\infty)\d x\}],
\end{align}
where the second equality also comes from \eqref{LK} except it is applied to $\Eb{e^{-\<{g,\theta}}}$ with the function $g(x):= - \log(\Eb{e^{-f(x+B_t-\lambda t}})$.

Note that for any measurable non-negative function $h$, by the change of variables $y=z+x-\lambda t$ and Fubini's theorem,
\begin{align*}
\int_{\R}\E[h(x+B_t-\lambda t)]e^{-2\lambda x}\d x=&\int_{\R}\int_{\R}h(x+z-\lambda t)\frac{1}{\sqrt{2\pi t}}e^{-\frac{z^2}{2t}}\d z e^{-2\lambda x}\d x\\
=&\int_\R \int_\R h(y)\frac{1}{\sqrt{2\pi t}}e^{-\frac{(y+\lambda t-x)^2}{2t}}\d y e^{-2\lambda x}\d x\\
=&\int_\R h(y)e^{-2\lambda y}\d y\int_\R \frac{1}{\sqrt{2\pi t}}e^{-\frac{(x+\lambda t-y)^2}{2t}}\d x\\
=&\int_\R h(y)e^{-2\lambda y}\d y. 
\end{align*}
Similarly, 
\[\int_{\R}\E[h(x+B_t-\lambda t)]\d x = \int_{\R}h(x)\d x.\] 
In particular, with $h(x)=1-e^{-f(x)}$,
\[
\int_{\R}\E[1-e^{-f(x+B_t-\lambda t)}](Z_\infty e^{-2\lambda x}+ Y_\infty)\d x =\int_{\R}[1-e^{-f(x)}](Z_\infty e^{-2\lambda x}+ Y_\infty)\d x
\]
which yields \eqref{lap-eq}.

\subsection{Characterization of the fixed points.}
We now turn to the main part of Theorem \ref{main-thm}, namely we prove that if $\theta$ is invariant, then there exist non-negative random variables $Z_\infty$ and $Y_\infty$ such that $\theta \overset{d}{=} \PPP((Z_\infty e^{-2\lambda x}+ Y_\infty)dx)$. 

We are going to show that for any $f\in C_c^+(\R)$,
\begin{equation}\label{goal}
\E[e^{-\<{f,\theta_t}}]=\E\left[\exp\left\{-\int_\R (1-e^{-f(y)})\left[e^{-2\lambda y}Z_\infty+Y_\infty\right]\d y\right\}\right]
\end{equation}
Recall that for any $t>0$, one has
\[
\E[e^{-\<{f,\theta_t}}]=\E\left[\exp\left\{\int_{-\infty}^\infty  \log(\E[e^{-f(x+B_t-\lambda t)}])\theta(\d x)\right\}\right].
\]
Let $K_f$ be such that the support of $f$ is contained in $[-K_f, K_f]$. Then
\begin{align*}
\E[1-e^{-f(x+B_t-\lambda t)}]\leq& \P(|x+B_t-\lambda t|\leq K_f)\nonumber\\
=&\int_\R \frac{e^{-\frac{z^2}{2t}}}{\sqrt{2\pi t}}\ind{|x+z-\lambda t|\le K_f}\d z\leq \frac{2K_f}{\sqrt{2\pi t}}= o_t(1)\,.
\end{align*}
As a consequence, $\E[e^{-\<{f,\theta_t}}]=\E[\exp(-(1+o_t(1))\Theta_t(f))]$ where 
\begin{equation}\label{exponent}
\Theta_t(f):=\int_{\R}\E[1-e^{-f(x+B_t-\lambda t)}]  \theta(\d x).
\end{equation}
From the structure of the known fixed points, we expect the main contribution to this integral to arise from points $x$ located either at $-\lambda t+o_t(1)$ (i.e. from the points corresponding to the expected $\PPP(Z_\infty e^{-2\lambda x} dx)$) or otherwise  located at $+\lambda t +o(1)$ (i.e. from the expected $\PPP(Y_\infty dx)$). We then split $\Theta_t(f)$ into the following three parts: 
\begin{equation}\label{exponent2}
\Theta_t(f)=\mathbf{Z}_t(f)+\mathbf{Y}_t(f)+E_t(f),
\end{equation}
where
\begin{align}
\mathbf{Z}_t(f):=&\int_{[-\lambda t-t^{2/3},-\lambda t+t^{2/3}]}\E[1-e^{-f(x+B_t-\lambda t)}]  \theta(\d x);\label{def-Z}\\
\mathbf{Y}_t(f):=&\int_{[\lambda t-t^{2/3},\lambda t+t^{2/3}]}\E[1-e^{-f(x+B_t-\lambda t)}]  \theta(\d x);\label{def-Y}\\
E_t(f):=&\int_{\R\setminus([-\lambda t-t^{2/3},-\lambda t+t^{2/3}]\cup [\lambda t-t^{2/3},\lambda t+t^{2/3}])}\E[1-e^{-f(x+B_t-\lambda t)}]  \theta(\d x).\label{def-E}
\end{align}
Observe that
\begin{align*}
&\E[1-e^{-f(x+B_t-\lambda t)}]\\
&=\E\left[(1-e^{-f(x+B_t-\lambda t)})\ind{x+B_t-\lambda t\ge - K_f}\right]\\
&=\E\left[(1-e^{-f(x+B_t-\lambda t)}) \vert x+B_t-\lambda t\ge - K_f\right]\P(x+B_t-\lambda t\ge - K_f).
\end{align*}
For any $a>0$, it is known that conditioned on $B_t\ge a_t=at+o(t)$,  $B_t-a_t$ converges in law to an exponential random variable with parameter $a$. One can even show a uniform convergence over $a_t\in [at-t^{2/3}, at+t^{2/3}]$. 

Let us recall the following basic estimate on the tail of normal distribution. For $z\to\infty$, one has
 \begin{equation}\label{keyBM}
 \P(B_1\geq z)=(1+o_z(1))\frac{1}{z\sqrt{2\pi}}e^{-z^2/2}.
 \end{equation}
It implies that for any fixed $K_f\ge0$, uniformly over $x\in[-\lambda t-t^{2/3}, -\lambda t+t^{2/3}]$, 
\begin{equation}\label{tailBM}
\P(x+B_t-\lambda t\geq -K_f)
=\frac{1+o_t(1)}{2\lambda\sqrt{2\pi t}}e^{-\frac{(\lambda t-x)^2}{2t}+2\lambda K_f}.
\end{equation}
Using again the notation $h(y)=1-e^{-f(y)}$. We see that for any $x\in[-\lambda t-t^{2/3}, -\lambda t+t^{2/3}]$, 
\begin{align*}
\E\left[(1-e^{-f(x+B_t-\lambda t)}) \ind{ x+B_t-\lambda t\ge - K_f}\right]&=\int_\R h(x+z-\lambda t)\frac{e^{-\frac{z^2}{2t}}}{\sqrt{2\pi t}}\ind{x+z-\lambda t\ge- K_f}\d z\\
&=\int_{-K_f}^\infty h(y) \frac{1}{\sqrt{2\pi t}}e^{-\frac{(y+\lambda t-x)^2}{2t}}\d y\\
&= \frac{1}{\sqrt{2\pi t}}e^{-\frac{(\lambda t-x)^2}{2t}}\int_{-K_f}^\infty h(y) e^{-\frac{y^2}{2t}-y (\lambda-\frac{x}{t})}\d y.
\end{align*}
By dominated convergence theorem, 
\[
\int_{-K_f}^\infty h(y) e^{-\frac{y^2}{2t}-y (\lambda-\frac{x}{t})}\d y\xrightarrow[t\to\infty]{}\int_{-K_f}^\infty h(y) e^{-2\lambda y}\d y
\]
as $t\to\infty$. Note that $h=1-e^{-f}$ vanishes on $(-\infty,-K_f)\cup (K_f,\infty)$ and that this convergence holds uniformly for $x\in[-\lambda t-t^{2/3}, -\lambda t+t^{2/3}]$. Moreover $\int_{-K_f}^\infty h(y) e^{-2\lambda y}\d y=\int_\R h(y) e^{-2\lambda y}\d y$. Going back to \eqref{def-Z} and \eqref{tailBM}, we get that
\begin{equation}\label{asymp-Z}
\mathbf{Z}_t(f)=(1+o_t(1))Z_t\int_\R (1-e^{-f(y)})   e^{-2\lambda y}\d y 
\end{equation}
where 
\begin{equation}\label{def-Z+}
Z_t:=\int_{-\lambda t-t^{2/3}}^{-\lambda t+t^{2/3}}\frac{1}{\sqrt{2\pi t}}e^{-\frac{(\lambda t-x)^2}{2t}}\theta(\d x)
\end{equation}
Note that by use of \eqref{tailBM}, we have
\begin{equation}\label{key-Z-p}
Z_t=(1+o_t(1))\frac{2\lambda}{e^{2\lambda}-1}\int_{-\lambda t-t^{2/3}}^{-\lambda t+t^{2/3}}\P(x+B_t-\lambda t\in [-1,0])\theta(\d x).
\end{equation}

\medskip
\ni
Now on the other side, for points $x\in [\lambda t-t^{2/3}, \lambda t+t^{2/3}]$, 
\begin{align*}
&\E[1-e^{-f(x+B_t-\lambda t)}]\\
&=\E\left[(1-e^{-f(x+B_t-\lambda t)})\ind{|x+B_t-\lambda t| \le K_f}\right]\\
&=\int_\R h(x+z-\lambda t)\frac{e^{-\frac{z^2}{2t}}}{\sqrt{2\pi t}}\ind{|x+z-\lambda t|\le K_f}\d z\\
&=\int_{-K_f}^{K_f} h(y) \frac{1}{\sqrt{2\pi t}}e^{-\frac{(y+\lambda t-x)^2}{2t}}\d y= \frac{1}{\sqrt{2\pi t}}e^{-\frac{(\lambda t-x)^2}{2t}}\int_{-K_f}^{K_f} h(y) e^{-\frac{y^2}{2t}-y (\lambda-\frac{x}{t})}\d y.
\end{align*}
Again, dominated convergence theorem shows that uniformly  for $x\in [\lambda t-t^{2/3}, \lambda t+t^{2/3}]$,
\[
\int_{-K_f}^{K_f} h(y) e^{-\frac{y^2}{2t}-y (\lambda-\frac{x}{t})}\d y\xrightarrow[t\to\infty]{}\int_{-K_f}^{K_f} (1-e^{-f(y)})\d y,
\]
where $\int_{-K_f}^{K_f} (1-e^{-f(y)})\d y=\int_\R (1-e^{-f(y)})\d y$ as $1-e^{-f}$ vanishes outside $[-K_f, K_f]$. We hence deduce that uniformly for $x\in [\lambda t-t^{2/3}, \lambda t+t^{2/3}]$, 
\begin{equation}\label{key-Y}
\E[1-e^{-f(x+B_t-\lambda t)}]=(1+o_t(1))\int_\R (1-e^{-f(y)})\d y \frac{1}{\sqrt{2\pi t}}e^{-\frac{(\lambda t-x)^2}{2t}}.
\end{equation}
Plugging it into \eqref{def-Y} yields that
\begin{equation}\label{asymp-Y}
\mathbf{Y}_t(f)= (1+o_t(1))Y_t\int_\R(1-e^{-f(y)})\d y
\end{equation}
where 
\begin{equation}\label{def-Z+}
Y_t:=\int_{\lambda t-t^{2/3}}^{\lambda t+t^{2/3}}\frac{1}{\sqrt{2\pi t}}e^{-\frac{(\lambda t-x)^2}{2t}}\theta(\d x).
\end{equation}
Notice that the previous arguments for \eqref{key-Y} also work if we replace $1-e^{-f}$ by $\ind{[-1, 0]}$.  It means that 
\begin{equation*} 
\P(x+B_t-\lambda t\in[-1,0])=(1+o_t(1))\frac{1}{\sqrt{2\pi t}}e^{-\frac{(\lambda t-x)^2}{2t}}
\end{equation*}
uniformly for $x\in [\lambda t-t^{2/3}, \lambda t+t^{2/3}]$. This implies that 
\begin{equation}\label{key-Y-p}
Y_t=(1+o_t(1))\int_{\lambda t-t^{2/3}}^{\lambda t+t^{2/3}}\P(x+B_t-\lambda t\in[-1,0])\theta(\d x)\,.
\end{equation}
\medskip
\ni
Next, for the term $E_t(f)$ in \eqref{def-E}, it is natural to split it into integrals on three disjoint intervals:
\begin{equation}\label{split-E}
E_t(f)=L_t(f)+C_t(t)+R_t(t)
\end{equation}
where $L_t(t):=\int_{-\infty}^{-\lambda t-t^{2/3}}\E[1-e^{x+B_t-\lambda t}]\theta(\d x)$, $C_t(t):=\int_{-\lambda t+t^{2/3}}^{\lambda t-t^{2/3}}\cdots$ and $R_t(f):=\int_{\lambda t+t^{2/3}}^\infty\cdots$. ($L_t,C_t,R_t$ respectively stand for Left, Center and Right terms). 

In view of \eqref{asymp-Z}, \eqref{asymp-Y}, \eqref{split-E}, \eqref{exponent2} becomes 
\begin{align*}
\Theta_t(f)&=(1+o_t(1))Z_t\int_\R (1-e^{-f(y)}) e^{-2\lambda y}\d y +(1+o_t(1)) Y_t\int_\R (1-e^{-f(y)})\d y\\
& \;\; +L_t(t)+C_t(t)+R_t(t).
\end{align*}
Here we claim the following assertions:
\begin{enumerate} \label{claim}
\item The random variables $Z_t$ are tight in $t>0$.
\item The random variables $Y_t$ are tight in $t>0$.
\item $L_t(f)$, $C_t(f)$, and $R_t(f)$ converge in probability to zero as $t \to \infty$.
\end{enumerate}
By admitting these three assertions, we deduce that as $t\to\infty$ along some subsequence, $\Theta_t(f)$ converges in law towards
\[
\int_\R(1-e^{-f(y)})\left[ Z_\infty  e^{-2\lambda y}+Y_\infty\right]\d y
\]
with some non-negative random variables $Z_\infty$ and $Y_\infty$ (who may not be independent). Recalling~\eqref{exponent}, this suffices to conclude that
\[
\E[e^{-\<{f,\theta}}]=\E\left[\exp\left\{-\int_\R (1-e^{-f(y)})\left[Z_\infty\, e^{-2\lambda y}+Y_\infty\right]\d y\right\}\right]
\]
which is what we want. 

It remains to verify the three assertions. We first state a basic concentration inequality for sum of independent Bernoulli random variables, which is follows from Chebyshev inequality and second Borel-Cantelli Lemma.
\begin{lemma}\label{concentration}
Let $(X_i, i\ge1)$ be a sequence of independent Bernoulli random variables such that $\E[X_i]=p_i\in(0,1]$. For any set $I\subset \N^*$, let $X_I:=\sum_{i\in I} X_i$. If $\E[X_I]<\infty$, then
\[
\P(|X_I-\E[X_I]|\ge \E[X_I]/2)\le \frac{4}{\E[X_I]}.
\]
If $\E[X_I]=\infty$, then $X_I=+\infty$ a.s.
\end{lemma}
Now we are ready to prove the above assertions (1),(2),(3).
Recall that we write $\theta=\sum_{i\in I}\delta_{x_i}$.

\subsection{The tightness of $Z_t$.} As mentioned in \eqref{key-Z-p}, we have 
\[
Z_t=(1+o_t(1))\frac{2\lambda}{e^{2\lambda}-1}\int_{-\lambda t-t^{2/3}}^{-\lambda t+t^{2/3}}\P(x+B_t-\lambda t\in [-1,0])\theta(\d x).
\]
The advantage of this expression is that if we define 
\[
\mathcal{Z}_t:=\sum_{i: x_i \in [-\lambda t-t^{2/3},-\lambda t+t^{2/3}]}\ind{x_i+B_t-\lambda t\in[-1,0]},
\]
then
\[
Z_t=(1+o_t(1)) \frac{2\lambda}{e^{2\lambda} -1}\E[\mathcal{Z}_t\vert \theta ].
\]
It thus suffices to prove the tightness of $\E[\mathcal{Z}_t\vert \theta ]$. Note that conditioned on $\theta$, $\mathcal{Z}_t$ is a sum of independent Bernoulli random variables. By use of Lemma \ref{concentration}, one sees that for any $K>0$,
\begin{align*}
\P\left(\E[\mathcal{Z}_t\vert \theta ]\ge K\right)
& \le  \P\left(\E[\mathcal{Z}_t\vert \theta ]\ge K; |\mathcal{Z}_t-\E[\mathcal{Z}_t\vert \theta ]|\ge \E[\mathcal{Z}_t\vert \theta ]/2\right)+\P\left(\mathcal{Z}_t\ge K/2\right)\\
& \le \frac{4}{K}+\P\left(\mathcal{Z}_t\ge K/2\right).
\end{align*}
Next, compare $\mathcal{Z}_t$ with $\theta_t$, we observe that $\mathcal{Z}_t\le \theta_t([-1,0])$. As $\theta_t([-1,0])$ has the  same law as $\theta([-1,0])$ which is finite a.s., 
\[
\P\left(\mathcal{Z}_t\ge K/2\right)\leq \P(\theta([-1,0])\ge K/2)
\]
It follows that
\[
\sup_{t>0}\P\left(\E[\mathcal{Z}_t\vert \theta ]\ge K\right)\le\frac{4}{K}+\P(\theta([-1,0])\ge K/2)\rightarrow 0,\textrm{ as } K\to\infty.
\]
This concludes the tightness of $Z_t$.

\subsection{The tightness of $Y_t$.}
In view of \eqref{key-Y-p}, the tightness of $Y_t$ follows from the same arguments as in the previous subsection.

\subsection{The convergences in probability of $L_t(f)$, $C_t(f)$, and $R_t(f)$.}
First, by \eqref{split-E}, one sees that
\begin{align}
L_t(f)\leq &\int_{-\infty}^{-\lambda t-t^{2/3}}\P(|x+B_t-\lambda t|\le K_f)\theta(\d x)=: L_t^+(f);\\
C_t(f)\le &\int_{-\lambda t+t^{2/3}}^{\lambda t -t^{2/3}}\P(|x+B_t-\lambda t|\le K_f)\theta(\d x)=: C_t^+(f)\\
R_t(f)\le &\int_{\lambda t+t^{2/3}}^{+\infty}\P(|x+B_t-\lambda t|\le K_f)\theta(\d x)=: R_t^+(f).
\end{align}
So, it suffices to show the convergence in probability towards zero for $L_t^+(f)$, $C_t^+(f)$, and $R_t^+(f)$.

First, let us introduce the random variables $L_s^+(f,t)$, $C_s^+(f,t)$ and $R_s^+(f,t)$ associated with some time $s>0$ (to be suitably chosen below) by replacing $\P(|x+B_t-\lambda t|\le K_f)\theta(\d x)$ by $\P(|x+B_s-\lambda s|\le K_f)\theta(\d x)$ in the integrals, as follows,
\[
L^+_s(f,t):=\int_{-\infty}^{-\lambda t-t^{2/3}}\P(|x+B_s-\lambda s|\le K_f)\theta(\d x)
\]

The following lemma  compares $L_t^+(f)$ and $L_s^+(f,t)$.
\begin{lemma}\label{ratio}
For $t$ large enough and $s= t + \sqrt{t}$ (respectively for $s= t-\sqrt{t}$), uniformly for $|x| \geq \lambda t + t^{2/3}$ (respectively for $|x|\leq \lambda t - t^{2/3}$), 
\[\frac{\P(|x+B_s-\lambda s|\leq K_f)}{\P(|x+B_t-\lambda t|\leq K_f)} \geq \frac12 e^{\frac{\lambda}{4} t^{1/6}}.\]
\end{lemma}
Its proof is postponed to the end. By means of this result, one sees that for $s=t+t^{1/2}$ and $t$ sufficiently large,
\[
L^+_s(t,f)\ge  \frac12 e^{\frac{\lambda}{4} t^{1/6}} L^+_t(f).
\]
So, for any $\eta>0$, one has
\begin{equation}\label{Lt-Ls}
\P(L_t^+(f)\ge \eta)\le \P\left(L^+_s(t,f)\ge  \frac12 e^{\frac{\lambda}{4} t^{1/6}} \eta\right).
\end{equation}
Note that if we set
\[
\calL^+_s(f,t):=\sum_{i\in I: x_i\in(-\infty, -\lambda t-t^{2/3})}\ind{|x_i+B_s^i-\lambda s|\le K_f}
\]
then $\E[\calL^+_s(f,t)\vert \theta]=L^+_s(f,t)$. Moreover, $\calL_s^+(f,t)$ is dominated by $\theta_s([-K_f,K_f])$ which is distributed as $\theta([-K_f,K_f])$. Applying Lemma \ref{concentration} to $\calL^+_s(f,t)$ conditioned on $\theta$ implies that for any $K>0$,
\begin{align}\label{key-bd-L}
&\P\left(L^+_s(t,f)\ge K\right)\nonumber\\
\le & \P\left(L^+_s(t,f)\ge K; |\calL^+_s(f,t)-L^+_s(f,t)|\ge \frac{L^+_s(f,t)}{2}\right)\nonumber\\
&\hspace{3cm}+\P\left(L^+_s(t,f)\ge K; \calL^+_s(f,t)\ge \frac{L^+_s(f,t)}{2}\right)\nonumber\\
\le & \frac{4}{K}+\P\left(\theta([-K_f,K_f])\ge K/2\right),
\end{align}
where the last inequality follows from the fact that $\calL^+_s(f,t)$ is stochastically dominated by $\theta([-K_f,K_f])$.

In view of \eqref{Lt-Ls} and \eqref{key-bd-L}, one gets that for any $\eta>0$,
\[
\P(L_t^+(f)\ge \eta)\le o_t(1)+\P\left(\theta([-K_f,K_f])\ge   e^{\frac{\lambda}{4}t^{1/6}}\eta/4\right)\xrightarrow[t\to\infty]{} 0,
\]
because $\theta$ is locally finite. This suffices to conclude that $L^+_t(f)$ converges in probability to zero. The same arguments can be applied to get the convergence in probability of $R^+_t(f)$.  For $C_t^+(f)$, we take $s=t-t^{1/2}$ instead of $t+t^{1/2}$, again the similar arguments hold and we obtain that $C_t^+(f)\xrightarrow{\P} 0$.

It remains to prove Lemma \ref{ratio}. Basic computation shows that
\begin{align*}
\P(|x+B_s-\lambda s|\le K_f)=&\int_\R \frac{1}{\sqrt{2\pi s}} e^{-\frac{z^2}{2s}}\ind{|x+z-\lambda s|\le K_f}\d z\\
=&\int_{-K_f}^{K_f} e^{-\frac{(y+\lambda s-x)^2}{2s}}\frac{\d y}{\sqrt{2\pi s}},
\end{align*}
and 
\[
\P(|x+B_t-\lambda t|\le K_f)=\int_{-K_f}^{K_f} e^{-\frac{(y+\lambda t-x)^2}{2t}}\frac{\d y}{\sqrt{2\pi t}}.
\]
Note that for $s=t+t^{1/2}$, for any $x$ such that $|x|\ge \lambda t+t^{2/3}$ and $y\in[-K_f, K_f]$, as long as $t$ is large enough, we have $s\le 2t$ and
\begin{align*}
-\frac{(y+\lambda s-x)^2}{2s}+\frac{(y+\lambda t-x)^2}{2t}=&\frac{(y-x)^2}{2st}t^{1/2}-\frac{\lambda^2}{2}t^{1/2}\\
\ge &\frac{\lambda t^{5/3+1/2}}{4t^2}=\frac{\lambda}{4} t^{1/6}.
\end{align*}
As a consequence, 
\[
\frac{\P(|x+B_s-\lambda s|\le K_f)}{\P(|x+B_t-\lambda t|\le K_f)}\ge \frac{1}{2} e^{\frac{\lambda}{4} t^{1/6}}\,.
\]
The case where $|x|\leq \lambda t - t^{2/3}$ and $s:=t-t^{1/2}$ is done similarly. This completes the proof.

\subsection*{Acknowledgements.}
The research of X.C is supported by ANR/FNS MALIN.
The research of C.G. and A.S.  is supported by the ERC grant LiKo 676999.

\bibliographystyle{alpha}
\bibliographystyle{acm}	
\bibliography{biblio}

\newcommand{\etalchar}[1]{$^{#1}$}
\begin{thebibliography}{AA{\etalchar{+}}09}

\bibitem[AA{\etalchar{+}}09]{AA09}
Louis-Pierre Arguin, Michael Aizenman, et~al.
\newblock On the structure of quasi-stationary competing particle systems.
\newblock {\em The Annals of Probability}, 37(3):1080--1113, 2009.

\bibitem[Bis17]{biskupBook}
Marek Biskup.
\newblock Extrema of the two-dimensional discrete gaussian free field.
\newblock In {\em PIMS-CRM Summer School in Probability}, pages 163--407.
  Springer, 2017.

\bibitem[BL16]{BL16}
Marek Biskup and Oren Louidor.
\newblock Extreme local extrema of two-dimensional discrete gaussian free
  field.
\newblock {\em Communications in Mathematical Physics}, 345, 2016.

\bibitem[Bov17]{bovier-book}
Anton Bovier.
\newblock {\em Gaussian Processes on Trees: From Spin Glasses to Branching
  Brownian Motion}.
\newblock Cambridge Studies in Advanced Mathematics. Cambridge University
  Press, 2017.

\bibitem[CD60]{ChoquetDeny}
G~Choquet and J~Deny.
\newblock Sur l'{\'e}quation de convolution $\mu\star\sigma=\mu$.
\newblock {\em CR Acad. Sci. Paris S{\'e}r. I Math}, 250:799--801, 1960.

\bibitem[CGS20]{CGS-BBM}
Xinxin Chen, Christophe Garban, and Atul Shekhar.
\newblock The fixed points of {B}ranching {B}rownian {M}otion.
\newblock {\em Preprint}, 2020.

\bibitem[Den60]{Deny}
Jacques Deny.
\newblock Sur l'{\'e}quation de convolution $\mu=\mu \star\sigma$.
\newblock {\em Seminaire Brelot-Choquet-Deny. Theorie du potentiel}, 4:1--11,
  1960.

\bibitem[Kal06]{kallenberg2006}
Olav Kallenberg.
\newblock {\em Foundations of modern probability}.
\newblock Springer Science \& Business Media, 2006.

\bibitem[Kal17]{kallenberg-book}
Olav Kallenberg.
\newblock {\em Random Measures, Theory and Applications}.
\newblock Probability Theory and Stochastic Modelling. Springer International
  Publishing, 2017.

\bibitem[Kin93]{kingman}
J.F.C. Kingman.
\newblock {\em Poisson Processes}.
\newblock Oxford science publications. Clarendon Press, 1993.

\bibitem[Lig78]{liggett78}
Thomas~M Liggett.
\newblock Random invariant measures for markov chains, and independent particle
  systems.
\newblock {\em Zeitschrift f{\"u}r Wahrscheinlichkeitstheorie und Verwandte
  Gebiete}, 45(4):297--313, 1978.

\bibitem[RA04]{AR05}
Anastasia Ruzmaikina and Michael Aizenman.
\newblock Characterization of invariant measures at the leading edge for
  competing particle systems.
\newblock {\em The Annals of Probability}, 33, 11 2004.

\bibitem[SZ17]{zeitouni2017}
Eliran Subag and Ofer Zeitouni.
\newblock The extremal process of critical points of the pure p-spin spherical
  spin glass model.
\newblock {\em Probability theory and related fields}, 168(3-4):773--820, 2017.

\end{thebibliography}

\end{document}